\documentclass[]{amsart} 
\usepackage[]{amsmath, amsthm, amsfonts}
\usepackage[all]{xy} 
\usepackage{hyperref}

\newcommand {\C} {{\mathbb C}}
  
\newcommand {\Z} {{\mathbb Z}}
 \newcommand {\Q} {{\mathbb Q}} 
\newcommand {\PP} {{\mathbb P}}
  
\newcommand {\F} {{\mathcal F}}

\newcommand {\dt} {{\bullet}}

\newcommand {\Y} {{\mathcal Y}}
\newcommand {\OO} {{\mathcal O}}

\newtheorem{thm}[subsection]{Theorem}
\newtheorem{cor}[subsection]{Corollary}
\newtheorem{lemma}[subsection]{Lemma}
\newtheorem{prop}[subsection]{Proposition}

\newtheorem{remark}[subsection]{Remark}
\newtheorem{ex}[subsection]{Example}

\begin{document}

 \title{Varieties with very little  transcendental cohomology}

 \author{ Donu Arapura} \thanks{Author partially supported by the NSF}
 \address{Department of Mathematics\\
   Purdue University\\
   West Lafayette, IN 47907\\
   U.S.A.}
 
 \maketitle

Given a complex smooth projective algebraic variety $X$,
 we define a natural number called the motivic
dimension $\mu(X)$ which is zero precisely when
all the cohomology of $X$ is generated by algebraic cycles. 
In general, it gives a measure of  the amount of 
transcendental cohomology of $X$. Alternatively,
 $\mu(X)$ may be rather loosely  thought of as measuring the complexity
of the motive of $X$, with Tate motives having $\mu=0$, motives
of curves having $\mu\le 1$ and so on.
Our interest in this notion stems from
the relation to the Hodge conjecture: it is easy to  see that it
holds for $X$ whenever $\mu(X)\le 3$. 
This paper contains a number  of  estimates  of
$\mu$; some elementary, some less so. 
With these estimates in hand, we conclude this paper by
checking or rechecking this conjecture in a number of examples:
uniruled fourfolds, rationally connected fivefolds, fourfolds
fibred by surfaces with $p_g=0$,  Hilbert schemes of a small number points on
surfaces with $p_g=0$, and generic hypersurfaces.


We will work over $\C$. Let $H^*(-)$
denote singular cohomology with rational coefficients.
The  motivic dimension  of a smooth projective variety
$X$ can be defined most succinctly in terms of
the length of the coniveau filtration on $H^*(X)$. We prefer to 
spell this out. The  {\em motivic dimension} $\mu(X)$ of $X$ is
the smallest integer $n$, such that  any
$\alpha\in H^i(X)$ vanishes on the complement of a Zariski closed
set all of whose components have codimension  at least $(i-n)/2$.
The meaning of this number is further clarified by the following:

\begin{lemma}
  Any $ \alpha\in H^i(X)$ can be decomposed as a finite sum of elements
of the form $f_{j*}(\beta_j)$, where $f_j:Y_j\to X $ are desingularizations of
subvarieties and  $\beta_j$ are classes of degree at most $\mu(X)$
on $Y_j$. In fact, $\mu(X)$ is the smallest integer such that the previous
statement holds for all $i\le \dim X$.
\end{lemma}

\begin{proof}
  The first statement follows from \cite[Cor. 8.2.8]{deligne}, and 
the second from this and the Hard Lefschetz theorem.
\end{proof}

\begin{cor}\label{cor:1}
   $\mu(X)=0$
if and only if all the cohomology of $X$ is generated by algebraic
cycles. A surface satisfies $\mu(X)\le 1$ if and only if $p_g(X)=0$.
We have
$$\\dim X\ge \mu(X)\ge \text{level}(H^*(X))=\max\{|p-q|\mid h^{pq}(X)\not=0\},$$
and the last inequality is equality if the generalized Hodge conjecture holds \cite{groth}.
\end{cor}

\begin{proof}
  The first statement is immediate. As for the second, the
Lefschetz $(1,1)$ theorem implies that $H^2(X)$ is spanned by
divisor classes if and only if $p_g=0$. The inequality $\mu(X)\ge \text{level}(H^*(X))$
follows from the fact that Gysin maps are morphisms of Hodge
structures up to Tate twists.
\end{proof}

My thanks to Su-Jeong Kang for her detailed comments. 

\section{Elementary Estimates}

 It will be convenient to extend the notion of motivic dimension to arbitrary complex algebraic
varieties, and for this we switch to homology.
Let $H_i(-)$ denote Borel-Moore homology with $\Q$ coefficients,
which is dual to compactly supported cohomology.  Translating
the above definition into homology leads to an integer $\mu^{big}(X)$,
 which we call the big motivic dimension that
makes sense for any variety $X$. So $\mu^{big}(X)$ is the smallest integer
such that every $\alpha\in H_i(X)$ 
lies in the image of some $f_*H_i(Y)\to H_i(X)$, where $Y$ is a
Zariski closed set whose components have dimension at most $
(i+\mu^{big}(X))/2$.
Of course, $\mu^{big}(X)$ coincides with $\mu(X)$ when $X$ is smooth and
proper, but it seems somewhat difficult to study in general.
It turns out to be more useful to restrict attention to certain cycles.
The identification  $H_i(X)\cong H^i_c(X)^*$ gives homology  a mixed Hodge
structure with weights $\ge -i$, i.e. $W_{-i-1}H_i(X)=0$
\cite{deligne}. We define  the
motivic dimension $\mu(X)$ by the replacing $H_i(X)$ by $W_{-i}H_i(X)=Gr_{-i}^WH_i(X)$
in the above definition of $\mu^{big}$. We have $\mu(X)\le \mu^{big}(X)$ with
equality when $X$ is smooth and proper.

\begin{prop}\label{prop:elem}
\begin{enumerate}
\item[]
\item[(a)] If $f:X'\to X$ is proper and surjective $\mu(X)\le\mu(X')$.

\item[(b)] If $Z\subset X$ is Zariski closed then
 $\mu(X)\le
\max(\mu(Z),\mu(X-Z))$
\item[(c)] If $\tilde X$ is a desingularization of a partial compactification
  $\bar X$ of
$X$, then $\mu(X)\le \mu(\bar X)\le \mu(\tilde X)$.
\item[(d)] $\mu(X_1\times X_2)\le \mu(X_1)+\mu(X_2)$.
\item[(e)] If $V\to X$ is a vector bundle then $ \mu(V)\le \mu(X)$
and $\mu(\PP(V))\le \mu(X)$.

\end{enumerate}
\end{prop}

\begin{proof}
By \cite[lemma 7.6, p. 110]{jannsen}, any element
 $\alpha\in W_{-i}H_{i}(X)$ lifts to an element of
 $\alpha'\in W_{-i}H_{i}(X')$. This in turn lies in $f_*
 W_{-i}H_{i}(Y')$ for some $f:Y'\to X'$ satisfying  $\dim Y'\le  (i+\mu(X'))/2$.
Therefore $\alpha$ lies in the image of $ W_{-i}H_{i}(f(Y'))$.
Since  $\dim f(Y')\le  (i+\mu(X'))/2$, therefore (a) holds.

Suppose $i\le m=\max(\mu(Z),\mu(X-Z))$.  We have an
 exact sequence of mixed Hodge structures 
 $$H_{i}(Z)\to H_{i}(X)\to H_{i}(X-Z)\to H_{i-1}(Z)$$
which can be deduced from \cite[prop. 8.3.9]{deligne}.
This implies by \cite[thm 2.3.5]{deligne} that
\begin{equation}
  \label{eq:W}
W_{-i}H_{i}(Z)\to W_{-i}H_{i}(X)\to W_{-i}H_{i}(X-Z)\to 0  
\end{equation}
is exact.
Given $\alpha\in W_{-i}H_i(X)$, let $\beta$ denote its image in $ W_{-i}H_i(X-Z)$.
Then $\beta = f_*(\gamma)$  for some $f:Y\hookrightarrow X-Z$
with $\dim Y\le  (i+m)/2$ and $\gamma\in W_{-i}H_i(Y)$.
 Let $\bar f:\bar Y\to X$ denote the closure
of $Y$. A sequence analogous to  \eqref{eq:W} shows that
$W_{-i}H_i(\bar Y)$ surjects onto $W_{-i}H_i(Y)$, therefore
$\gamma$ extends to a class $\bar\gamma \in
W_{-i}H_i(\bar Y)$. The difference
$\alpha-\bar f_*(\bar \gamma)$ lies in the image of $W_{-i}H_{i}(Z)$,
and therefore in $g_*H_i(T)$ for some $g:T\to Z$ with  $\dim T\le
(i+m)/2$. This proves (b).

Let $\bar X$ be a partial compactification of $X$. Then as above,
we see that any class in $W_{-i}H_{i}(X)$ extends to $\bar X$.
Therefore $\mu(X)\le \mu(\bar X)$. The remaining inequality of (c)
follows from (a).

Statement (d) follows from the K\"unneth formula
$$W_{-i}H_i(X_1\times X_2) = \bigoplus_{j+k=i} W_{-j}H_j(X_1)\otimes W_{-k}H_k(X_2)$$

Finally for (e), let $r=rk(V)$. Suppose that the Gysin images of
$W_{-i+2r}(Y_j)$ span $W_{-i+2r}H_{i-2r}(X)$ and satisfy
$\dim Y_j\le (\mu(X)+i-2r)/2$.  Then $W_{-i}H_i(V|_{Y_j})$ will span $W_{-i}H_{i}(V)$ by the
  Thom isomorphism theorem. The inequality $\mu(\PP(V))\le \mu(X)$ can
  be proved by a similar argument. We omit the details
since we will show something more general in corollary~\ref{cor:brauer}.
\end{proof}

\begin{cor}
  If $X=\cup X_i$ is given as a finite disjoint union of locally
  closed subsets, $\mu(X)\le \max \{\mu(X_i)\}$. In particular,
$\mu(X)=0$ if $X$ is a disjoint union of open subvarieties of  affine spaces.
\end{cor}

The last statement, which is of course well known, implies that flag varieties and
toric varieties have $\mu=0$.

\begin{cor}
If $X_1$ and $X_2$ are birationally equivalent smooth projective
varieties, then $\mu(X_2)\le \max(\mu(X_1),\dim X_1-2)$
\end{cor}

\begin{proof}
Since an iterated blow up of  $X_1$ dominates $X_2$,
it is enough to prove this when $X_2$ is the blow up of $X_1$
along a smooth centre $Z$. Then
$$\mu(X_2) \le \max(\mu(X_1-Z), \mu(\PP(N))) \le
\max(\mu(X_1),\mu(Z))$$
where $N$ is the normal bundle of $Z$.
\end{proof}

Recall that a variety is uniruled if it has a rational curve passing
through the general point.

\begin{cor}\label{cor:uni}
  If $X$ is uniruled, then $\mu(X)\le \dim X-1$.
\end{cor}

\begin{proof}
 By standard arguments \cite{kollar}, $X$ is dominated by 
a blow up of $Y\times \PP^1$, where $\dim Y=\dim X-1$. 
\end{proof}

\begin{cor}
If $X$ is a smooth projective variety with a $\C^*$-action,
$\mu(X)$ is less than or equal to the dimension of the
fixed point set.
\end{cor}

\begin{proof}
Bialynicki-Birula \cite{bb} has shown that $X$ can be decomposed
into a disjoint union of vector bundles over components of the fixed
point set. The corollary now follows from the previous results.
\end{proof}

From this, one recovers the well known fact that the Hodge numbers $h^{pq}(X)$ vanish
when $|p-q|$ exceeds the dimension of the fixed set.

\begin{prop}\label{prop:bloch}
Suppose that $X$ is a smooth projective variety such that 
the Chow group of  zero cycles $CH_0(X)\cong \Z$. Then 
 $\mu(X)\le \dim X-2$.
\end{prop}

\begin{proof}
This follows from the theorem of Bloch-Srinivas \cite{bs}
that a positive multiple of the diagonal $\Delta\subset X\times X$
 is rationally equivalent to a sum $\xi\times X+ \Gamma$,
 where $\xi\in Z_0(X)$ is a zero cycle,  
 and $\Gamma$ is supported on $X\times D$
 for some divisor  $D\subset X$.  This implies
 that the identity map on cohomology $H^i(X)$ factors through
 the Gysin map $H^{i-2}(\tilde D)\to H^i(X)$ for $i>0$ and a 
 resolution of singularities $\tilde D\to D$.
\end{proof}

Recall that a projective variety is rationally connected if any two general points
can be connected by a rational curve. Examples include hypersurfaces
in $\PP^n$ with degree less than $n+1$, and more generally Fano varieties
\cite{kollar}.

\begin{cor}
If $X$ is rationally connected then $\mu(X)\le \dim X-2$.
\end{cor}

\begin{proof}
Rational connectedness forces $CH_0(X)\cong \Z$.
\end{proof}

\section{Estimates for Fibrations}

\begin{thm}\label{thm:fibre}
  Suppose that $f:X\to S$ is a smooth projective morphism.
 Then 
$$\mu(X)\le \max_{s\in S(\C)}{\mu(X_s)}+ \dim S$$
\end{thm}

\begin{proof}
We prove this  by induction on $d=\dim S$. 
Let $m= \max_{s\in S(\C)}{\mu(X_s)}$. Let  $H_1, H_2,\ldots$ denote
irreducible components of the relative
Hilbert scheme  $Hilb_{X/S}$ which surject onto $S$.
Choose a desingularization $\tilde
\F_k\to \F_{k,red}$ of the reduced universal
family over each $H_k$.
$$
\xymatrix{
 \tilde \F_k\ar[r]\ar[rdd] & \F_k\ar[d]\ar@{^{(}->}[r] & X\times_S Hilb_{X/S}\ar[d] \\ 
  & H_k\ar[d]\ar@{^{(}->}[r] & Hilb_{X/S}\ar[ld] \\ 
  & S & 
}
$$
Let $disc(\tilde\F_k\to S)\subseteq S$ denote discriminant i.e. the complement of the maximal open set
over which this maps is smooth.
Then $T= \bigcup_k disc(\tilde \F_k\to S)\subset S$ is a countable union of proper
subvarieties, therefore its complement is  nonempty by Baire's theorem.
Choose $s\in S-T$. By assumption, there exists
a finite collection of subvarieties $Y_{ij}\subset X_s$ such that
$\dim Y_{ij}\le (m+i)/2$ and their images generate $H_i(X_s)$.
For each $Y_{ij}$, we choose one of the families $\F_k$ containing
it and rename it $\Y_{ij}$; likewise set $H_{ij}=H_k$ and $\tilde \Y_{ij}=\tilde\F_k$.

Choose an open set $U\subset S-\bigcup disc(\Y_{ij}\to S)$ containing $s$.
Therefore $X|_U\to U$ and each $\tilde \Y_{ij}|_U\to U$ are smooth and
thus topological fibre
bundles. Consequently the images of $\tilde \Y_{ij,t}$ will generate $H_i(X_t)$ for
any $t\in U$. After replacing $H_{ij}$ by the
normalization of $S$ in $H_{ij}$, and $\mathcal{Y}_{ij}$ by the fibre product,
we can assume that $H_{ij}\to S$ is finite over its image.
After shrinking $U$ if necessary and replacing all maps by their
restrictions to $U$, we can assume that  the maps
$g_{ij}:\tilde{\mathcal{Y}}_{ij}\to U$ are  still smooth, 
and that $U$ is nonsingular and affine. The last assumption implies
that for any local system $L$ of $\Q$-vector spaces,
we have
\begin{equation}
  \label{eq:HiL}
  H_c^i(U,L)= H^{2d-i}(U, L^*)^* = 0
\end{equation}
for $i<d$ since $U$ is homotopic to a CW complex of dimension at most
$d$ \cite[thm 1.22]{voisin}.
The Gysin
images of $\tilde{\mathcal{Y}}_{ij,t}$ generate the homology of $X_t$
for each $t\in U$, or dually the cohomology of $H^i(X_t)$ injects
into $\oplus_j H^i(\tilde{\mathcal{Y}}_{ij,t})$. Since the monodromy
actions are semisimple \cite[thm 4.2.6]{deligne}, the map of
local systems $R^if_*\Q|_U\to \oplus_j R^ig_{i,j,*}\Q$ is split
injective.
Thus 
\begin{equation}
  \label{eq:HckU}
H_c^k(U,R^if_*\Q)\to \bigoplus_j H_c^k(U,R^ig_{i,j,*}\Q)  
\end{equation}
is injective. Since the Leray spectral sequence degenerates
\cite{deligne-L},
we get an injection 
$$H_c^p(f^{-1}U)\to \bigoplus_{i\le p-d, j}
H_c^p(\tilde{\mathcal{Y}}_{ij})$$
Note that the bound $i\le p-d$ follows from (\ref{eq:HiL}).
We therefore have a surjection 
 $$ \bigoplus_{i\le p-d, j} H_p(\tilde{\mathcal{Y}}_{ij})\to H_p(f^{-1}U) $$
 Since 
 \begin{eqnarray}
   \label{eq:Yij}
\dim \tilde{\mathcal{Y}}_{ij}\le \frac{m+i}{2}+d\le \frac{m+d+p}{2}   
 \end{eqnarray}
 for $i\le p-d$, we have $\mu(f^{-1}(U))\le \mu^{big}(f^{-1}U)\le
 m+d$. By induction $\mu(f^{-1}(S-U))\le m+d$. Thus
 $\mu(X)\le \max(\mu(f^{-1}U),\mu(f^{-1}(S-U)))\le m+d$ as required.
 \end{proof}

 \begin{cor}
 If  a smooth projective variety $X$ can be covered by a family of surfaces whose general
member is smooth with $p_g=0$. Then $\mu(X)\le \dim X-1$.
\end{cor}

\begin{proof}
By assumption, there is a family $Y\to T$ of surfaces with $p_g=0$ and a dominant
map $\pi: Y\to X$. After restricting to a subfamily,
we can assume that $\pi$ is generically  finite.
So that $\dim T= \dim X-2$, and therefore
$\mu(X)\le \dim X-1$ by
corollary~\ref{cor:1} and theorem~\ref{thm:fibre}. 
\end{proof}

It is worth noting that many standard examples of surfaces with
$p_g=0$ lie in families. So there are
nontrivial examples of varieties admitting  fibrations of the above
type. We give an explicit class of examples generalizing   Enriques
surfaces (cf. \cite[p. 184]{bpv}).

\begin{ex}
   Fix $n\ge 2$.
  Let $i:(\PP^1)^n\to (\PP^1)^n$ be the involution which acts by
  $[x_0,x_1]\mapsto [-x_0,x_1]$ on each factor. Choose a divisor
  $B\subset (\PP^1)^n$ defined by a general $i$-invariant polynomial of
  multidegree $(4,4,\ldots 4)$, and let $\pi:X\to (\PP^1)^n$ be the double cover
  branched along $B$. The involution $i$ can be seen to lift to a fixed point free
  involution of $X$. Let $S$ be the quotient. $X$ is covered by a
  family of  an $K3$ surfaces $\pi^{-1}((\PP^1)^2\times {t})$ which
  induces a family of Enriques surfaces on $S$. Thus $\mu(S)\le n-1$.
  On the other hand  $X$ can be checked to be Calabi-Yau, and thus the Kodaira
  dimension of $S$ equals $0$. Therefore it cannot be uniruled. So this
  estimate on $\mu(S)$ does not appear to follow from the previous bounds.
\end{ex}

In view of proposition~\ref{prop:elem} (d),
we may hope for a stronger estimate
$$\mu(X)\le \max_{s\in S(\C)}{\mu(X_s)}+ \mu(S)$$
Unfortunately it may fail without extra assumptions: 

\begin{ex}
Let $S$ be an Enriques surface which can be
realized as the quotient of a K3 surface $\tilde S$ by a fixed point free
involution $\sigma$. Let $\sigma$ act on $\PP^1\times \PP^1$
by interchanging factors. Define $X=(\PP^1\times \PP^1\times \tilde S)/\sigma$.
The natural map $X\to S$ is an etale locally trivial $\PP^1\times \PP^1$-bundle.
An easy calculation shows that $level(H^2(X)) =2$, while $p_g(S)=q(S)=0$.
Thus $\mu(X)=2 > \mu(\PP^1\times \PP^1)+\mu(S)=0$.
\end{ex}

 \begin{thm}\label{thm:fibre2}
 Suppose that $f:X\to S$ is a smooth projective morphism over a
 quasiprojective base.
   If the monodromy action of $\pi_1(S)$ on $H^*(X_s)$ is trivial, then
$$\mu(X)\le \max_{s\in S(\C)}{\mu(X_s)}+ \mu(S)$$
 \end{thm}

Before proving this, we make some general remarks.
If $f:X\to S$ is a not necessarily smooth projective morphism, then
there is a filtration  $L^\dt H^i(X)\subset H^i(X)$ called the Leray filtration
associated to the Leray spectral sequence. This is a filtration by
sub mixed Hodge structures  \cite[cor. 4.4]{arapura}. When $f$ is also
smooth the spectral sequence degenerates \cite{deligne-L}, so that
$Gr^p_L H^{p+q}(X)\cong H^p(S, R^qf_*\Q)$ carries a mixed Hodge structure.

\begin{lemma}\label{lemma:leray}
Suppose that $f:X\to S$ is a smooth projective map, and let  $m= \max_{s\in S(\C)}{\mu(X_s)}$.
 Then   $Gr^W_{k+i}H_c^k(S, R^if_*\Q)$ injects into  $Gr^W_{k+i} Gr^k_LH_c^{k+i}(Y)$ for some
Zariski closed  $Y\subset X$ satisfying $\dim Y\le
\frac{m+i}{2}+\dim S$.
\end{lemma}

 \begin{proof}
 As in the proof of theorem~\ref{thm:fibre}, we can find a nonempty
 affine open set $U\subset S$ and a morphism of smooth $U$-schemes $\tilde \Y=\cup_j \tilde \Y_{ij}\to f^{-1}U$
whose fibres generated the homology of $H_i(X_s)$. The lemma holds
when $S$ is replaced by $U$ with $Y=im\tilde \Y$ thanks to 
 (\ref{eq:HckU}) and (\ref{eq:Yij}). Let $\bar\Y\subset X$ denote the
 closure of $im \Y$.  By induction, we have a subset
 $\Y'\subset f^{-1}(S-U)$ which satisfies the lemma over $S-U$.
We have a commutative diagram with exact rows
$$
\xymatrix{
   &H_c^{k+i}(f^{-1}U)\ar[r]\ar[d] & H_c^{k+i}(X)\ar[r]\ar[d] & H_c^{k+i}(f^{-1}(S-U))\ar[d] \\ 
 0\ar[r]& H_c^{k+i}(\Y)\ar^{j}[r] & H_c^{k+i}(\Y'\coprod \bar\Y)\ar[r] & H_c^{k+i}(\Y'\coprod (\bar\Y-\Y))
}
$$
Applying $Gr^W_{k+i}Gr_L^k$ and making appropriate identifications
results in a commutative diagram 
$$
\xymatrix{
   Gr^W_{k+i}H_c^k(U, R^if_*\Q)\ar[r]\ar^{\iota}[d] & Gr^W_{k+i}H_c^k(S, R^if_*\Q)\ar[r]\ar^{\iota''}[d] & Gr^W_{k+i}H_c^k(S-U, R^if_*\Q)\ar^{\iota'}[d] \\ 
  Gr^W_{k+i}H_c^{k+i}(\Y)\ar^{j}[r] & Gr^W_{k+i}H_c^{k+i}(\Y'\coprod\bar \Y)\ar[r] & Gr^W_{k+i}H_c^{k+i}(\Y'\coprod (\bar\Y-\Y))
}
$$
The top row is  exact, but the
bottom row need not be. Nevertheless
$j$ is  injective since $Gr^W_{k+i}Gr_L^k$ preserves
injections. The maps $\iota$ and $\iota'$ are also injective by the above discussion.
The injectivity of $\iota''$ follows a diagram chase. Thus $Y= \Y'\cup
\bar\Y$ does the job.
 \end{proof}

\begin{proof} 
Let $m= \max_{s\in S(\C)}{\mu(X_s)}$.
The sheaves $R^if_*\Q$ are constant, so we have an isomorphism
 \begin{equation}\label{eq:L}
H_c^k(S,R^if_*\Q)\cong H_c^k(S)\otimes H^i(X_s)
\end{equation}
as vector spaces.  As already noted,
the Leray spectral sequence for $f$ degenerates yielding a mixed Hodge
structure on the left side.
We claim that  \eqref{eq:L}  can be made compatible with  mixed Hodge
structures, at least after taking the associated graded with respect to the
weight filtration. We have a surjective morphism of pure polarizable Hodge structures
$Gr^W_*H_c^i(X)\to H^i(X_s)$, which admits a right inverse $\sigma$
since this category is semisimple \cite[lemma 4.2.3]{deligne}.
 The image of $H_c^k(S)$ lies
in $L^kH^k_c(X)$. Thus $Gr^WGr_L(f^*\otimes \sigma)$ induces the desired
identification
\begin{equation}
  \label{eq:Gr}
Gr^W_* H_c^k(S)\otimes  H^i(X_s)\cong Gr^W_*H_c^k(X,R^if_*\Q)
\end{equation}
Moreover, this is canonical in the sense that it  is compatible with
base change with respect to any morphism $T\to S$.

 Let $S_k\to S$ be
 a map such that $\dim S_k\le (\mu(S)+k)/2$ and such that
 $W_{-k}H_k(S_k)$ generates $W_{-k}H_k(S)$ (note that $S_k$ may have
 several components).
Dually $Gr^W_{k}H_c^k(S)$ injects into $Gr^W_{k}H_c^k(S_k)$. 
Thus there are injections
$$
\xymatrix{
  Gr^W_{k}H_c^k(S)\otimes H^i(X_s)\ar@{^{(}->}[r]\ar[d]^{\cong} &  Gr^W_{k}H_c^k(S_k)\otimes H^i(X_s)\ar[d]^{\cong} \\ 
 Gr^W_{k+i}H^k_c(S, R^if_*\Q)\ar@{^{(}->}[r] & Gr^W_{k+i}H^k_c(S_k, R^if_*\Q)
}
$$
By  lemma~\ref{lemma:leray} $ Gr^W_{k+i}H^k_c(S_k, R^if_*\Q)$ 
 injects into some $ Gr^W_{i+k}Gr_L^kH_c^{k+i}(\Y_{ki})$ with 
$$\dim \Y_{ki}\le \dim S_k +\frac{m+i}{2}\le \frac{m+\mu(S)+ k+i}{2}$$
Combining this  with the degeneration of
Leray, shows that  the map
$$ Gr^W_pH_c^{p}( X)\to \bigoplus_{k+i= p} Gr^W_pH_c^p(\Y_{ki})$$
is injective, and hence we have a surjection
$$  \bigoplus_{k+i= p} Gr^W_{-p}H_p(\Y_{ki})\to
Gr^W_{-p}H_p(X)$$

\end{proof}

\begin{cor}\label{cor:zariski}
With the above notation, suppose that there exists a nonempty
Zariski open set $U\subseteq S$ such that 
 $X|_U\to U$ is topologically a product. Then $\mu(X)\le \mu(X_s)+\mu(S)$.
\end{cor}

\begin{proof}
$\pi_1(U)\to \pi_1(S)$ is surjective.
\end{proof}

\begin{cor}\label{cor:brauer}
  If $X\to S$ is a Brauer-Severi morphism (i.e. a smooth map whose fibres  are projective spaces)
then $\mu(X)\le\mu(S)$.
\end{cor}

\begin{proof}
  Since $\dim H^i(\PP^N)\le 1$, the  monodromy representation is
trivial. So $\mu(X)\le \mu(\PP^N)+\mu(S)$.
\end{proof}
\section{Symmetric powers }

In this section, we give our take on Abel-Jacobi theory. When $X$ is a
smooth projective curve, the symmetric powers $S^nX$ are projective
bundles over the Jacobian $J(X)$ for $n\gg 0$. This implies that the motivic
dimension of $S^nX$ stays bounded. We consider what happens for more
general smooth projective varieties. We note that $S^nX$ are singular in
general, but only mildly so. These are in the class of $V$-manifolds or orbifolds,
which  satisfy Poincar\'e duality
with rational coefficients, hard  Lefschetz and purity of  mixed Hodge structures.
So  for our purposes, we can treat them as smooth.
In particular,  we work with the original cohomological definition of $\mu$.
When $X$ is a surface, the Hilbert schemes provides natural desingularization of the
symmetric powers, and we give estimates for these as well.

 We can identify $H^*(S^nX)$ with the $S_n$-invariants of $H^*(X^n)$. Let 
 $$sym:H^*(X^n)\to H^*(X^n)^{S_n}\cong H^*(S^nX)$$
 denote the symmetrizing operator $\frac{1}{n!}\sum \sigma$. 

\begin{lemma}\label{lemma:alphabeta}
Let $f:X\to Y$ be a morphism of smooth projective varieties.
If every class in $H^*(X)$ is of the form $(f^*\alpha)\cup \beta$ where $\beta$ is an
algebraic cycle, then $\mu(X)\le \mu(Y)$.
\end{lemma}

\begin{proof}
Let $N^pH^i(X)$ denote the span of the Gysin images of desingularizations of
subvarieties of codimension $\ge p$.  Then it is enough to show that $N^pH^i(X)=H^i(X)$
for some $p\ge (i-\mu(Y))/2$. Any element of $H^i(X)$ can be written as $f^*\alpha\cup \beta$
where $\beta\in H^{2k}(X)$ is an algebraic cycle and $\alpha\in H^{i-2k}(Y)$.
This implies that $\alpha\in N^qH^{i-2k}(Y)$ and $\beta\in N^kH^{2k}(X)$ for some  $q$ satisfying $q\ge (i-2k-\mu(Y))/2$. 
By \cite{ak2}, $f^*\alpha\cup \beta\in N^{q+k}H^i(X)$. Therefore $p=q+k$ gives the desired value.
\end{proof}

\begin{cor}\label{cor:alphabeta}
Suppose that $G$ is a finite group.
Let $f:X\to Y$ be an equivariant  morphism of smooth projective varieties with
$G$-actions satisfying the above assumption, then $\mu(X/G)\le \mu(Y/G)$.
\end{cor}

\begin{proof}
This is really a corollary of the proof which proceeds as above
with the identifications $N^pH^i(X/G) = N^pH^i(X)^G$ etcetera.
\end{proof}

\begin{cor}\label{cor:alphabeta2}
Let $f:X\to Y$ be a morphism of smooth projective varieties, which is a Zariski
locally trivial fibre bundlle with fibre $F$ satisfying $\mu(F)=0$. Then
$\mu(S^nX)\le \mu(S^nY)$ for all $n$.
\end{cor}

\begin{proof}
The conditions  imply that the hypotheses  of the previous corollary holds
for $X\to Y$ as well as its symmetric powers. 
\end{proof}

\begin{thm}
Let $X$ be a smooth projective variety.
\begin{enumerate}
\item[(a)] $\mu(S^nX)\le n \mu(X)$.
\item[(b)] If $\mu(X) \le 1$, then the sequence $\mu(S^nX)$ is bounded.
If $\dim X\le 2$ then this is bounded above by $h^{10}(X)$.
\item[(c)] If $level(H^{2*}(X))>0$, then $level(H^*(S^nX))$ and 
$\mu(S^nX)$ are unbounded.
\end{enumerate}
\end{thm}

\begin{proof}
Inequality (a) follows from proposition \ref{prop:elem}.

Suppose that $\mu(X)\le 1$. Then $H^*(X)$ is spanned by algebraic cycles
and classes $f_{j*}\beta$
with $f_j:Y_j\to X$ and $\beta\in H^{1}(Y_j)$.  Let $Y$ be the disjoint union
of $Y_j$, and let $f^n:S^nY\to S^nX$ denote the natural map.
Fix a basis $\beta_1\ldots \beta_{N}$ of $H^1(Y)$.
 Then  $H^*(S^nX)$ is
 spanned by classes of the form 
 \begin{equation}\label{eq:beta}
 [f^n_*sym(p_m^*(\beta_{i_1}\times\ldots \times \beta_{i_m}))]\cup \gamma
\end{equation}
 where $p_m:Y^n\to Y^m$ is a projection onto the first $m$ factors, and $\gamma$ is an
 algebraic cycle.  Notice that the expression in \eqref{eq:beta} vanishes if any of the $\beta_{i_j}$'s are repeated. Therefore we can assume that $m\le N$. 
 Let $g_j$ denote the inclusions of components of  the algebraic cycle $f^{n*}\gamma$.
 We can rewrite the expressions in
 \eqref{eq:beta} as 
 $$f^n_*(sym(\ldots)\cup f^{n*}\gamma)\in \text{span}\{f^n_*g_{j_*}(g_j^*sym(\ldots))\}$$
 This shows that $H^*(S^nX)$ is spanned by  Gysin images of classes of degree at most $N$.
 Thus $\mu(S^nX)\le N$.  This proves the first part of (b).
 
 When $X$ is a curve of genus $g=h^{10}(X)$, the Abel-Jacobi map $S^nX\to J(X)$
 can be decomposed into a union of  projective space bundles. Proposition  \ref{prop:elem} implies
  that $\mu(S^nX)\le g=\dim J(X)$.  
  
  Suppose $X$ is a surface, then $\mu(X)\le 1$ forces
  $p_g(X)=0$. If $q=h^{10}(X)=0$ then $\mu(X)=0$, so $\mu(S^nX)=0$ by (a).
So we may assume that $q>0$, which implies that the Albanese map $\alpha:X\to Alb(X)$
is nontrivial. If $\dim \alpha(X)=2$, it is easy to see that the pullback of a generic two form
from $Alb(X)$ would be nontrivial. This would imply  that the image is a curve. Let $C$ be
the normalization of $\alpha(X)$, then we get a map $\phi:X\to C$. 
Note the the genus $g$ of $C$ is necessarily equal to $q$ since
$\alpha^*$ induces an isomorphism on the space $1$-forms and it
factors as $H^0(\Omega_{Alb(X)}^1)\to H^0(\Omega_C^1)\hookrightarrow H^0(\Omega_X^1)$.
With the help of
the hard Lefschetz theorem, we can see that every cohomology class on $X$ is
of the form $(\phi^*\beta)\cup \gamma$ where $\gamma$ is an algebraic cycle.
Likewise for $S^nX$. Therefore by corollary~\ref{cor:alphabeta} and previous
paragraph $\mu(S^nX)\le \mu(S^nC) \le g$.

Suppose that $\alpha\in H^{ij}(X)$ is a nonzero class with $i\not= j$ and $i+j$ even. Then
$sym(\alpha^{\otimes n})$ provides a nonzero class in $H^{ni,nj}(S^nX)$, which shows
that the level and hence the motivic dimension go to infinity.
\end{proof}

We give an example where $\mu(X)>1$ but $\mu(S^nX)$ stays bounded.

\begin{ex}\label{ex:cy}
Let $X$ be a rigid Calabi-Yau threefold. (A number of
such  examples are known, cf.  \cite{schoen}.)
Then  the  Hodge numbers satisfy $h^{10}=h^{20}=h^{21}=0$ and $h^{30}=1$. 
It follows that $H^2(X)$ and $H^4(X)$ are generated by algebraic cycles.
Choose a generator $\alpha\in H^{30}(X)$ normalized so that the class 
$\delta=sym(\alpha\times \bar\alpha)\in H^6(S^2X)$ is rational. The space 
of $S_2$-invariant classes of type $(3,3)$ in $H^3(X)\otimes H^3(X)$ is one dimensional.
Therefore $\delta$ must coincide with the K\"unneth component of the diagonal $\Delta$ in 
$H^3(X)\otimes H^3(X)$. The remaining K\"unneth components of $\Delta$ are necessarily algebraic.
Therefore $\delta$ is algebraic. 
Thus  all factors 
$$sym(H^3(X)^{\otimes k}\otimes H^{2i_1}(X)\otimes\ldots \otimes H^{2i_{n-k}}(X))$$
 in $H^i(S^nX)$ are spanned by algebraic cycles if $k$ is even, or
 $\{sym((\xi\times (\text{alg. cycle})\mid \xi\in H^3(X)\}$ if $k$ is odd.
It follows, by a modification of lemma \ref{lemma:alphabeta}, that $\mu(S^nX)\le 3$
for all $n$.
\end{ex}

 We review the basic facts about Hilbert schemes of surfaces. Proofs and references 
can be found in the first 30 or so pages of \cite{gott}.
When $X$ is a smooth projective surface, there is a natural desingularization
$\pi_n:Hilb^nX\to S^nX$ given by the (reduced) Hilbert scheme of $0$-dimensional
subschemes of length $n$. 
The symmetric product $S^nX$ can be decomposed into a disjoint union of locally closed sets
$$\Delta_{(\lambda_1,\ldots \lambda_k)} = \{\lambda_1x_1+\ldots \lambda_kx_k \mid x_i\not= x_j\}$$
indexed by partitions $\lambda_1\ge \lambda_2\ldots$ of $n$, where the
elements of $S^nX$ are written additively.
Let $\tilde \Delta_\lambda\subset X^k$ denote the preimage under
$(x_1,\ldots, x_k)\mapsto \lambda_1x_1+\ldots \lambda_kx_k$.
 The map $\tilde \Delta_\lambda\to \Delta_\lambda$
is etale with Galois group $G$. The  group can be described explicitly
 by grouping the terms in the partition as follows:
$$\underbrace{\lambda_1=\ldots =\lambda_{d_1}}_{d_1}
>\underbrace{\lambda_{d_1+1}=\ldots =\lambda_{d_1+d_2}}_{d_2}>\ldots \, 
\underbrace{\ldots\lambda_{d_1+d_2+\ldots d_\ell}}_{d_\ell}=\lambda_k>0$$
Then $G=S_{d_1}\times \ldots S_{d_\ell}$ and
$ \Delta_\lambda$ is isomorphic to  an open subset $S^{d_1}X\times \ldots S^{d_\ell}X$.
We record the following key fact \cite[lemma 2.1.4]{gott}

\begin{lemma}
$\pi_n^{-1}\Delta_\lambda\times_{\Delta_\lambda} \tilde \Delta_\lambda\to \tilde \Delta_\lambda$
is isomorphic to a Zariski open subset of the $k$-fold product $\prod H_{\lambda_i}\to X^k$.
Each $H_{\lambda_i}\to X$
is a Zariski locally trivial fibre bundle whose fibre can be identified
with the subscheme  $Hilb_0^{\lambda_i}\C^2\subset Hilb^{\lambda_i}\C^2$
parameterizing schemes supported at the origin of $\C^2$.
\end{lemma}


We note that $Hilb^{\lambda_i}_0\C^2$ is smooth and
$\C^*$ acts on it with isolated fixed points.
Therefore the fibres have $\mu=0$. 
We can see that $G$ acts  on $\pi_n^{-1}\Delta_\lambda\times_{\Delta_\lambda} \tilde \Delta_\lambda$
by permuting the $H_{\lambda_i}$'s. Consequently 
\begin{equation}\label{eq:Delta}
\pi_n^{-1}\overline{\Delta}_\lambda\cong S^{d_1}H_{\lambda_{d_1}}\times  S^{d_2}H_{\lambda_{d_1+d_2}}\times\ldots S^{d_\ell}H_{\lambda_{d_1+\ldots d_\ell}}
\end{equation}

\begin{prop}
 If $X$ is a surface with $p_g=0$ then $\mu(Hilb^nX) \le \min(n,\sqrt{2n}q)$ for all $n$,
 where $q=h^{10}$. In particular, $\mu(Hilb^nX) = 0$ when $q=0$.
\end{prop}
 
\begin{proof}
  The closures of each of the strata $\Delta_\lambda$ are dominated by $X^k$,
so that  $\mu( \Delta_\lambda )\le k\le n$. The maps $\pi^{-1}\Delta_\lambda\to \Delta_\lambda$ are
Zariski locally trivial fibrations with fibres having $\mu=0$. Thus $\mu(Hilb^nX) \le n$ follows from this together with proposition \ref{prop:elem} and corollary \ref{cor:zariski}.

Each $H_{\lambda_i}\to X$ is a bundle  with $\mu=0$ fibres.
 Therefore $\mu(S^{d_i}H_{\lambda_{d_1+\ldots d_i}}) \le q$ by
 corollary \ref{cor:alphabeta2} and the previous theorem.
Combing this with \eqref{eq:Delta}, we see that
 $\Delta_\lambda$  has motivic dimension at most $\ell q$.
To estimate $\ell$, we use
\begin{eqnarray*}
n &=&d_1 \lambda_{d_1}+ d_{2} \lambda_{d_1+d_{2}}+\ldots d_\ell \lambda_{d_1+\ldots d_\ell}\\
&\ge& \lambda_{d_1}+  \lambda_{d_1+d_{2}}+\ldots \lambda_{d_1+\ldots d_\ell}\\
&\ge& \ell +(\ell-1)+\ldots 1
\end{eqnarray*}
to obtain $\ell\le \sqrt{2n}$.  This gives the remaining inequality 
$\mu(Hilb^nX) \le \sqrt{2n}q$ 

\end{proof}

\section{Applications to the Hodge conjecture}

Jannsen \cite{jannsen} has extended the Hodge conjecture to an
arbitrary variety $X$. This states that
 any class in $Hom(\Q(i), W_{-2i}H_{2i}(X))$, which should be thought of as a Hodge
 cycle, is a linear combination of fundamental classes of $i$-dimensional subvarieties.
Lewis \cite{lewis} has given a similar extension for the generalized Hodge conjecture
which would say that an irreducible sub Hodge structure of
$W_{-i}H_i(X)$ with
level at most $\ell$ should lie in the Gysin
image of a subvariety of dimension bounded by $(\ell + i)/2$. 
Both statements are equivalent to the usual
forms for smooth projective varieties.

\begin{prop}

  \begin{enumerate}
  \item[]
\item[(a)] If $\mu(X)\le 3$, then the Hodge conjecture holds for $X$.
\item[(b)] If $\mu(X)\le 2$, then the generalized Hodge conjecture
holds for $X$.
  \end{enumerate}
  
\end{prop}

\begin{proof}
Suppose that $\mu(X)\le 3$.
Then  any Hodge class $\alpha\in H_{2i}(X)$ lies in $f_* W_{-2i}H_{2i}(Y)$  for some subvariety 
with components satisfying $\dim Y_j\le (3+2i)/2$.
Let $\tilde Y=\cup \tilde Y_j$ be a  desingularization of a compactification of $Y$.
An argument similar to the proof of proposition \ref{prop:elem} shows
that the natural map 
$$W_{-2i}H_{2i}(\tilde Y)\to f_* W_{-2i}H_{2i}(Y)$$
 is a surjective
morphism of polarizable Hodge structures. This map admits a section, since
the category of such structures is semisimple \cite{deligne}.
It follows that $\alpha$ can be lifted to a Hodge cycle $\beta$ on $\tilde Y$.
This can be viewed as a Hodge cycle in cohomology under the Poincar\'e
duality isomorphism $H_{2i}(\tilde Y )= \oplus_jH^{2\dim Y_j-2i}(\tilde Y_j)(\dim Y_j)$.
Since  $2\dim Y_j-2i\le 3$, this forces the degree of $\beta$ to be $0$ or $2$.
consequently $\beta$ must be an algebraic cycle. Hence the same is true
for its image $\alpha$.

The second statement  is similar. 
With notation as above, a sub Hodge structure of $H_i(X)$ 
is the image of a sub Hodge structure of
 $\oplus_jH^{2\dim Y_j-i}(\tilde Y_j)(\dim Y_j)$ with  $2\dim
 Y_j-i\le 2$.
Since the generalized Hodge conjecture is trivially true in this range,
this structure is contained in the Gysin image of map from a
subvariety of expected dimension.
\end{proof}

\begin{cor}[Conte-Murre]
  The Hodge (respectively generalized Hodge) conjecture holds
 for uniruled fourfolds (respectively threefolds).
\end{cor}

\begin{cor}[Laterveer]
  The Hodge (respectively generalized Hodge) conjecture
holds for rationally connected smooth projective fivefolds (respectively fourfolds).
\end{cor}

\begin{cor}
If $X$ is a smooth projective variety with a  $\C^*$-action, then
the Hodge (respectively generalized Hodge) conjecture holds
if the fixed point set has dimension at most $3$ (respsectively $2$).
\end{cor}

The last result can also be deduced from the main theorem of
\cite{ak}, which says in effect that these conjectures factor through
 the Grothendieck group of varieties. 
The class of $X$ in the Grothendieck group can be expressed as a linear combination of classes
of components of the fixed point set times Lefschetz classes.

\begin{cor}
  The Hodge (respectively generalized Hodge) conjecture holds
for a smooth projective fourfold (respectively threefold)
which can be covered by a family of surfaces whose general
member is smooth with $p_g=0$.
\end{cor}

\begin{cor}
Let $X$ be a smooth projective surface with $p_g=0$.  
Then the  Hodge (respectively generalized Hodge) conjecture
holds for $S^nX$ for any $n$ if $q\le 3$ (respectively $q\le 2$).
The Hodge  conjecture
holds for $Hilb^nX$ if  $n \le 3$, or if $n\le 7$ and $q=1$, or if $q=0$.
The generalized Hodge conjecture holds for $Hilb^nX$ if
$n= 2$, or if $n\le 3$ and $q=1$, or if $q=0$.
\end{cor}

As a final example, suppose that
$X$ is a rigid Calabi-Yau variety. Then we saw that $\mu(S^nX)\le 3$ in example \ref{ex:cy}.
So the Hodge conjecture holds for $S^nX$. (We are being a little circular in our logic, 
since we essentially verified the conjecture in the course of estimating $\mu(S^nX)$.)

We have a Noether-Lefschetz result in this setting.
For the statement, we take  ``sufficiently general'' to mean
that the set of exceptions forms a countable union of proper Zariski
closed subsets of the parameter space.

\begin{thm}
  Let $X\subset \PP^N$ be a smooth projective variety such that
the Hodge (respectively generalized Hodge) conjecture holds. Then
 there exists an effective constant $d_0>0$ such that the Hodge (respectively generalized
 Hodge)  conjecture holds
for $H\subset X$ when $H\in \PP(\OO_{X}(d))$ is a sufficiently general
hypersurface  and $d\ge d_0$.

\end{thm}

\begin{remark}
The effectivity of $d_0$ depends on having enough information about
$X\subset\PP^N$. As will be clear from the proof, it would be
sufficient to know the Chern
classes of $X, \OO_X(1)$, the Castelnuovo-Mumford regularity of
$\OO_X$ and $T_X$, and $h^0(T_X)$.
 For $X=\PP^N$, we have $d_0=N+1$.
\end{remark}

\begin{proof}
Let $n=\dim X-1$. We $d_0$ to be the smallest integer for which
\begin{enumerate}
\item[(a)] $d_0\ge 2$.
\item[(b)]  $ h^{n0}(H)>h^{n0}(X)$
for nonsingular $H\in \PP(\OO_{X}(d))$ with $d\ge d_0$. 
(This ensures that $H^n(H)/im H^n(X)$ has length $n$; a fact needed
at the end.)
\item[(c)] $H\in \PP(\OO_{X}(d))$, $d\ge d_0$, has nontrivial moduli, at least infinitesimally
\end{enumerate}
To clarify these conditions, note that from the exact sequence
$$0\to \omega_X\to \omega_X(d)\to \omega_{H}\to 0$$
and Kodaira's vanishing theorem we obtain
\begin{eqnarray*}
  h^0(H, \omega_{H}) &= &h^0(X,
  \omega_X(d))-h^0(\omega_X)+h^1(\omega_X)\\
 &= &\chi(\omega_X(d))-h^0(\omega_X)+h^1(\omega_X)\\
\end{eqnarray*}
The right side is explicitly computable by Riemann-Roch, so we get an
effective lower bound for (b). For (c), we require
$H^1(T_{H})\not=0$. Once we choose $d_0$ so that $H^1(T_X(-d))=0$ for $d\ge
d_0$, standard exact sequences yield
\begin{eqnarray*}
h^1(T_{H}) &\ge& h^0(\OO_{H}(d))-h^0(T_X|_H)  \\
&\ge& h^0(\OO_{X}(d))-h^0(T_X)-1
\end{eqnarray*}
We can compute the threshold which makes the right side positive.
When $X=\PP^N$, we  see by a direct computation that
assumptions (a),(b), (c) are satisfied as soon as $d\ge N+1$.

Let $U_d\subset \PP(\OO(d))$ denote  the
set of nonsingular hypersurfaces in  $X$.
 For $H\in U_d$,
the weak Lefschetz theorem guarantees an isomorphism $H^i(X)\cong
H^i(H)$ for $i<n$, so the (generalized) Hodge conjecture
holds for these groups by assumption.
This together with the hard Lefschetz theorem takes care of the
range $i>n$.
So it remains to treat $i=n$.
We have an orthogonal decomposition
\begin{equation}
  \label{eq:van}
  H^n(H) = im\, H^n(X)\oplus V_H
\end{equation}
as Hodge structures, where $V_H$ is the kernel of the Gysin map $H^n(H)\to
H^{n+2}(X)$ \cite[sect. 2.3.3]{voisin}. Equivalently, after
restricting to a Lefschetz pencil in $\PP(\OO(d))$ containing $H$,
 $V_H$ is just the space of vanishing cycles [loc. cit.].
As $H$ varies,
$H^n(H)$ determines a local system over $U_d$, and  (\ref{eq:van}) 
is compatible with the monodromy action.  The action of $\pi_1(U_d)$ 
respects the intersection pairing $\langle,\rangle$. So the image
of a finite index subgroup $\Gamma\subset \pi_1(U_d)$ lies in the
identity component $ Aut(V_H,\langle,\rangle)^o$, which is
a symplectic or special orthogonal group according to the parity of $n$. By
 \cite[thm 4.4.1]{deligne-W2}, the image of $\Gamma$ 
 is either finite or  Zariski dense in $ Aut(V_H,\langle,\rangle)^o$. The first
possibility can be ruled out by checking that the Griffiths period
map on $U_d$ is nontrivial
\cite[lemma 3.3]{ct}. The nontriviality of the period map is
guaranteed by Green's Torelli theorem \cite[thm
0.1]{green} and our assumption (c).
Thus the Zariski closure
$\overline{im(\pi_1(U_d))}^{Zar}\supseteq Aut(V_H,\langle,\rangle)^o$.
To finish the argument, we recall that the
Mumford-Tate group $MT(V_H)\subset GL(V_H)$ is an algebraic subgroup which
 leaves all sub Hodge structures invariant. By  \cite[prop
7.5]{deligne-K3}, when $H$ is sufficiently general
$MT(V_H)$ must contain a finite index subgroup of the
image of $\pi_1(U_d)$, and therefore $Aut(V_H,\langle,\rangle)^o$.  Thus 
$MT(V_H)$ contains the symplectic or special orthogonal group.
In either case $V_H$ is an irreducible Hodge structure of length $n$. The
(generalized) Hodge conjecture concerns  sub Hodge
structures of $H^n(H)$ of length less than $n$. So
these must come from $X$, and therefore be contained in images
 of Gysin maps of subvarieties of expected codimension by our initial
assumptions.
\end{proof}

\begin{remark}
For  the Hodge conjecture alone,  the bound on $d_0$
can be  improved. 
 When $n$ is odd, we may take $d_0=1$ since the
Hodge conjecture is vacuous for $H^n$.
When $n$ is even, we can choose the smallest $d_0$ so that assumptions (a)
and (c) hold, since the mondromy action on
$V_H$ is irreducible and nontrivial (by the Picard-Lefschetz formula) and therefore
it cannot contain a Hodge cycle.
\end{remark}

\begin{cor}
   Let $X\subset \PP^N$ be a smooth projective variety.
  \begin{enumerate}
  \item[(a)] If $\mu(X)\le 3$, then the Hodge conjecture holds
for  every sufficiently general
hypersurface section of degree $d\gg 0$. This is in particular the
case, when $X$ is a rationally connected fivefold.
 \item[(b)] If $\mu(X)\le 2$, then the generalized Hodge conjecture holds
for  every sufficiently general
hypersurface section of degree $d\gg 0$. 
  \end{enumerate}
\end{cor}

\end{document}